\documentclass{article}
\usepackage{amsfonts}

\input{tcilatex}
\begin{document}

\begin{center}
\textbf{Some Properties of Mannheim Curves in Galilean and Pseudo - Galilean
space}

\textbf{Alper Osman \"{O}\u{G}RENM\.{I}\c{S}, Handan \"{O}ZTEK\.{I}N, Mahmut
ERG\"{U}T}

F\i rat University, Science Faculty, Mathematics Department

23119 Elaz\i \u{g} / T\"{U}RK\.{I}YE

ogrenmisalper@gmail.com, handanoztekin@gmail.com, mergut@firat.edu.tr

\bigskip

\textbf{Abstract}
\end{center}

\begin{quote}
The aim of this work is to study the Mannheim curves in \linebreak
3-dimensional Galilean and Pseudo - Galilean space. We obtain the
characterizations between the curvatures and torsions of the Mannheim
partner curves.
\end{quote}

\linebreak Mannheim curve, Curvature, Torsion.

\textit{2000 AMS Classification:}\textbf{\ }53A35; 53B30.

\bigskip

\begin{center}
\textbf{1.} \textbf{Introduction. }
\end{center}

In classical differential geometry, there are a lot of studies on Bertrand
and Mannheim curves. Many interesting results on Bertrand and Mannheim
curves have been obtained by many mathematicians ( see [2,7,8,9,10,11,13] ).
We can see in most studies, a characteristic property of Bertrand and
Mannheim curves which asserts the existence of a lineer relation between
curvature and torsions. There are a lot of works on Bertrand curves but
there are rather a few works on Mannheim curves.

In this study, we have done a study about some special curves in Galilean
and Pseudo - Galilean space. However, to the best of author's knowledge,
Mannheim curves has not been presented Galilean and Pseudo - Galilean space
in depth. Thus, the study is proposed to serve such a need.

Our study is organized as follows. In secion 2, some fundamental \linebreak
properties of Galilean and Pseudo - Galilean space are given which will be
used in the later sections. In section 3, we give some characterizations of
Mannheim curves in Galilean space. We also give some properties of Mannheim
curves in Pseudo - Galilean space in section 4.

\begin{center}
\textbf{2. Basic notions and properties}
\end{center}

The Galilean space is a three dimensional complex projective space $P_{3}$
in which the absolute figure $\{w,f,I_{1},I_{2}\}$ consists of a real plane $%
w$ (the absolute plane), a real line $f\subset w$ (the absolute line) and
two complex conjugate points $I_{1},I_{2}\in f$ (the absolute points).

We shall take, as a real model of the space $G_{3}$ , a real projective
space $P_{3}$ with the absolute $\{w,f\}$ consisting of a real plane $%
w\subset G_{3}$ and a real line $f\subset w$ on which an elliptic involution 
$\varepsilon $ has been defined.We introduce homogeneous coordinates in $%
G_{3}$ in such a way that the absolute plane $w$ is given by $x_{0}=0,$ the
absolute line $f$ by $x_{0}=x_{1}=0$ and elliptic involution by $%
(x_{0}:x_{1}:x_{2}:x_{3})=(1:x:y:z),$ the distance between the points $%
P_{i}=(x_{i},y_{i},z_{i}),$ $i=1,2,$ is defined by

\[
d(P_{1},P_{2})=\{%
\begin{array}{c}
\qquad \quad \qquad \qquad \left\vert x_{2}-x_{1}\right\vert ,\qquad
if\qquad x_{1}\neq x_{2} \\ 
\sqrt{(y_{2}-y_{1})^{2}+(z_{2}-z_{1})^{2}}\qquad if\qquad x_{1}=x_{2}%
\end{array}%
. 
\]

Let $\alpha $ be a curve in $G_{3},$ defined by arclength $\alpha
:I\rightarrow G_{3}$ and parametrized by the invariant parameter $s\in I$ ,
given in the coordinate form%
\[
\alpha (s)=(s,\text{ }y(s),\text{ }z(s)). 
\]

Then the curvature $\kappa (s)$ and the torsion $\tau (s)$ are defined by%
\[
\kappa (s)=\sqrt{y^{\prime \prime ^{2}}(s)+z^{\prime \prime ^{2}}(s)} 
\]%
\[
\tau (s)=\frac{\det (\alpha ^{\prime }(s),\text{ }\alpha ^{\prime \prime
}(s),\text{ }\alpha ^{\prime \prime \prime }(s))}{\kappa ^{2}(s)} 
\]%
and associated moving trihedron is given by%
\begin{eqnarray*}
T(s) &=&\alpha ^{\prime }(s)=(1,\text{ }y\prime (s),\text{ }z^{\prime }(s))
\\
N(s) &=&\frac{1}{\kappa (s)}\alpha ^{\prime \prime }(s)=\frac{1}{\kappa (s)}%
(0,\text{ }y^{\prime \prime }(s),\text{ }z^{\prime \prime }(s)) \\
B(s) &=&\frac{1}{\kappa (s)}(0,-z^{\prime \prime }(s),\text{ }y^{\prime
\prime }(s)).
\end{eqnarray*}

The vectors $T,$ $N,$ $B$ are called the vectors of the tangent, principal
normal and binormal line of $\alpha ,$ respectively. For their derivatives
the following Frenet formulas hold%
\begin{eqnarray}
T^{\prime }(s) &=&\kappa (s)N(s)  \nonumber \\
N^{\prime }(s) &=&\tau (s)B(s)  \TCItag{2.1} \\
B^{\prime }(s) &=&-\tau (s)N(s).  \nonumber
\end{eqnarray}

More about the Galilean geometry can be foun in [12].

The Pseudo-Galilean \ geometry is one of the real Cayley-Klein
geometries\linebreak (of projective signature $(0,0,+,-)$, explained in [3].
The absolute of the Pseudo-Galilean \ geometry is an ordered triple $%
\{w,f,I\}$ where \ $w$ is the ideal (absolute) plane, $f$ is line in $w$ and 
$I$ is the fixed hyperbolic involution of points of $f.$

As in [3], Pseudo-Galilean scalar product can be written as

\[
\left\langle v_{1},v_{2}\right\rangle =\left\{ 
\begin{array}{c}
x_{1}x_{2}\qquad \quad \text{if }x_{1}\neq 0\vee x_{2}\neq 0 \\ 
y_{1}y_{2}-z_{1}z_{2}\text{ if }x_{1}=0\wedge x_{2}=0%
\end{array}%
\right. 
\]%
where $v_{1}=(x_{1},y_{1},z_{1}),$ $v_{2}=(x_{2},y_{2},z_{2})$.

It leaves invariant the Pseudo-Galilean norm of the vector $v=(x,y,z)$
defined by 
\[
\left\Vert v\right\Vert =\left\{ 
\begin{array}{c}
x,\text{ }x\neq 0 \\ 
\sqrt{\left\vert y^{2}-z^{2}\right\vert },\text{ }x=0.%
\end{array}%
\right. 
\]

A vector $v=(x,y,z)$ is said to be non-isotropic if \linebreak $x\neq 0.$
All unit non-isotropic vectors are the form $(1,y,z)$. \ For \ isotropic
\linebreak vectors $x=0$ holds. There are four types of isotropic vectors:
spacelike \linebreak $(y^{2}-z^{2}>0),$ timelike $(y^{2}-z^{2}<0)$ and two
types of lightlike $(y=\pm z)$ \linebreak vectors. A non-lightlike isotropic
vector is a unit vector if $y^{2}-z^{2}=\pm 1.$

Let $\alpha :I\longrightarrow G_{3}^{1},$ $I\subset \mathbb{R}$ be a curve
given by 
\[
\alpha (t)=(x(t),y(t),z(t)), 
\]%
where $x(t),y(t),z(t)$ $\in C^{3}$ (the set of three times continuosly
differentiable functions) and $t$ run through a real interval [3].

A curve $\alpha (t)$ in $G_{3}^{1}$ is said to be admissible curve if $%
x^{\prime }(t)\neq 0$ [3].

\bigskip

The curves in $G_{3}^{1}$ are characterized as follows [4,5]

\bigskip

\textbf{Type I:} Let $\alpha $ be an admissible curve in $G_{3}^{1},$
parametrized by arclength $t=s,$ given in coordinate form%
\[
\alpha (s)=(s,\text{ }y(s),\text{ }z(s)). 
\]

Then the curvature $\kappa _{\alpha }(s)$ and the torsion $\tau _{\alpha
}(s) $ are defined by%
\[
\kappa _{\alpha }(s)=\sqrt{\left\vert y^{\prime \prime ^{2}}(s)-z^{\prime
\prime ^{2}}(s)\right\vert } 
\]%
\[
\tau _{\alpha }(s)=\frac{\det (\alpha ^{\prime }(s),\text{ }\alpha ^{\prime
\prime }(s),\text{ }\alpha ^{\prime \prime \prime }(s))}{\kappa _{\alpha
}^{2}} 
\]%
and associated moving trihedron is given by%
\begin{eqnarray*}
T_{\alpha }(s) &=&\alpha ^{\prime }(s)=(1,\text{ }y\prime (s),\text{ }%
z^{\prime }(s)) \\
N_{\alpha }(s) &=&\frac{1}{\kappa _{\alpha }(s)}\alpha ^{\prime \prime }(s)=%
\frac{1}{\kappa _{\alpha }(s)}(0,\text{ }y^{\prime \prime }(s),\text{ }%
z^{\prime \prime }(s)) \\
B_{\alpha }(s) &=&\frac{1}{\kappa _{\alpha }(s)}(0,\text{ }z^{\prime \prime
}(s),\text{ }y^{\prime \prime }(s)).
\end{eqnarray*}

The vectors $T_{\alpha },$ $N_{\alpha },$ $B_{\alpha }$ are called the
vectors of the tangent, \linebreak principal normal and binormal line of $%
\alpha ,$ respectively. For their derivatives the \linebreak following
Frenet formulas hold%
\begin{eqnarray}
T_{\alpha }^{\prime }(s) &=&\kappa _{\alpha }(s)N_{\alpha }(s)  \nonumber \\
N_{\alpha }^{\prime }(s) &=&\tau _{\alpha }(s)B_{\alpha }(s)  \TCItag{2.2} \\
B_{\alpha }^{\prime }(s) &=&\tau _{\alpha }(s)N_{\alpha }(s)  \nonumber
\end{eqnarray}

\textbf{Type II:} Let $\beta $ be an admissible curve in $G_{3}^{1},$
parametrized by arclength $s,$ given in coordinate form%
\[
\beta (s)=(s,y(s),0). 
\]

Then the curvature $\kappa _{\beta }(s)$ and the torsion $\tau _{\beta }(s)$
are defined by%
\begin{eqnarray*}
\kappa _{\beta }(s) &=&y^{\prime \prime }(s) \\
\tau _{\beta }(s) &=&\frac{a_{2}^{\prime }(s)}{a_{3}(s)}
\end{eqnarray*}%
and associated moving trihedron is given by%
\begin{eqnarray*}
T_{\beta }(s) &=&(1,\text{ }y\prime (s),\text{ }0), \\
N_{\beta }(s) &=&(0,\text{ }a_{2}(s),\text{ }a_{3}(s)), \\
B_{\beta }(s) &=&(0,\text{ }a_{3}(s),\text{ }a_{2}(s)).
\end{eqnarray*}%
where, $y$, $a_{2},$ $a_{3}\in C^{\infty },$ $s\in I\subseteq \mathbb{R}.$

The vectors $T_{\beta },$ $N_{\beta },$ $B_{\beta }$ are called the vectors
of the tangent, principal normal and binormal line of $\beta ,$
respectively. For their derivatives the following Frenet's formulas hold%
\begin{eqnarray}
T_{\beta }^{\prime }(s) &=&\kappa _{\beta }(s)(\cosh \phi (s)N_{\beta
}(s)-\sinh \phi (s)B_{\beta }(s)),  \nonumber \\
N_{\beta }^{\prime }(s) &=&\tau _{\beta }(s)B_{\beta }(s),  \TCItag{2.3} \\
B_{\beta }^{\prime }(s) &=&\tau _{\beta }(s)N_{\beta }(s),  \nonumber
\end{eqnarray}%
where $\phi (s)$ is the angle between $a(s)=(0,a_{2}(s),$ $a_{3}(s))$ and
the plane $z=0.$

\begin{center}
\textbf{3. Mannheim Curves in Galilean Space}
\end{center}

\textbf{Definition 3.1. }Let $\alpha $ and $\alpha _{1}$ be curves in
3-dimensional Galilean space $G_{3}.$ If there exists a corresponding
relationship between the space curves $\alpha $ and $\alpha _{1}$ such that,
at the corresponding points of the curves, the principal normal lines of $%
\alpha $ coincides with the binormal lines of $\alpha _{1}$, then $\alpha $
is called a Mannheim curve and $\alpha _{1}$ is called a Mannheim partner
curve of $\alpha $. The pair $\{\alpha ,\alpha _{1}\}$ is said to be a
Mannheim pair. [6]

\textbf{Remark 3.2. }Let $\alpha $ be a curve in Galilean space $G_{3}.$
Then $\alpha $ is Mannheim curve if and only if its curvature $\kappa
_{\alpha }$ and torsion $\tau _{\alpha }$ satisfy the relation \linebreak $%
\kappa _{\alpha }=c\tau _{\alpha }^{2}$ for some constant $c.$ [10]

\textbf{Theorem 3.3.} Let $\alpha $ be a Mannheim curve in Galilean space $%
G_{3}.$ Then $\alpha _{1}$ is the Mannheim partner curve of $\alpha .$ Then
the curvature $\kappa _{1}$ and the torsion $\tau _{1}$ of $\alpha _{1}$
satisfy the following equation

\begin{equation}
\tau _{1}^{^{\prime }}=\frac{\kappa _{1}}{\lambda }(\lambda ^{2}\tau
_{1}^{2}+1)  \tag{3.1}
\end{equation}%
for some nonzero constant $\lambda .$

Proof. Suppose that $\alpha (s)$ is a Mannheim curve in Galilean space $%
G_{3}.$ Then we can write

\begin{equation}
\alpha (s_{1})=\alpha _{1}(s_{1})+\lambda (s_{1})B_{1}(s_{1})  \tag{3.2}
\end{equation}%
for some function $\lambda (s_{1}).$ By taking derivative of (3.2) with
respect to $s_{1}$ and using the Frenet equations in Galilean space $G_{3},$
we get

\begin{equation}
T\frac{ds}{ds_{1}}=T_{1}+\lambda ^{^{\prime }}B_{1}-\lambda \tau _{1}N_{1}. 
\tag{3.3}
\end{equation}%
Since $B_{1}$ is coincident with $N$, we obtain

\[
\lambda ^{^{\prime }}(s_{1})=0 
\]%
that means that $\lambda $ is nonzero constant. Thus, we have

\begin{equation}
T\frac{ds}{ds_{1}}=T_{1}-\lambda \tau _{1}N_{1}.  \tag{3.4}
\end{equation}%
On the other hand we have

\begin{equation}
T=T_{1}\cos \theta +N_{1}\sin \theta  \tag{3.5}
\end{equation}%
where $\theta $ is the angle between $T$ and $T_{1}$ at the corresponding
points of $\alpha $ and $\alpha _{1}.$ Differentiating of (3.5) with respect
to $s_{1}$, we get

\begin{eqnarray}
\kappa N\frac{ds}{ds_{1}} &=&-T_{1}\theta ^{^{\prime }}\sin \theta +\kappa
_{1}N_{1}\cos \theta +N_{1}\theta ^{^{\prime }}\cos \theta +\tau
_{1}B_{1}\sin \theta  \nonumber \\
\kappa N\frac{ds}{ds_{1}} &=&-T_{1}\theta ^{^{\prime }}\sin \theta
+N_{1}(\kappa _{1}+\theta ^{^{\prime }})\cos \theta +\tau _{1}B_{1}\sin
\theta .  \TCItag{3.6}
\end{eqnarray}%
Since $\{\alpha ,\alpha _{1}\}$ is a Mannheim pair, we obtain

\[
\kappa _{1}+\theta ^{^{\prime }}=0 
\]%
and therefore we have

\begin{equation}
\theta ^{^{\prime }}=-\kappa _{1}.  \tag{3.7}
\end{equation}%
From

\bigskip (3.4) and (3.5), we find that

\begin{equation}
\lambda \tau _{1}=-\tan \theta .  \tag{3.8}
\end{equation}%
By taking derivative of this last equation and applying (3.7), we obtain

\begin{equation}
\lambda \tau _{1}^{^{\prime }}=-\theta ^{^{\prime }}(1+\tan ^{2}\theta ). 
\tag{3.9}
\end{equation}%
If we consider (3.7) and (3.8) in (3.9), we get

\[
\tau _{1}^{^{\prime }}=\frac{\kappa _{1}}{\lambda }(\lambda ^{2}\tau
_{1}^{2}+1). 
\]%
Hence the proof is completed.

\textbf{Proposition 3.4.} Let $\alpha $ be a Mannheim curve in Galilean
space $G_{3}$ and $\alpha _{1}$ be the Mannheim partner curve of $\alpha .$
If $\alpha $ is a generalized helix, then $\alpha _{1}$ is a planar curve.

Proof. Let $T,$ $N,$ $B$ the tangent , principal normal and binormal vector
field of the curve $\alpha $, respectively. From the properties of
generalized helix and the definition of Mannheim curves in Galilean space $%
G_{3},$ we have

\[
\kappa <N,P>=0 
\]%
and 
\[
\kappa <B_{1},P>=0 
\]%
for a constant direction $P$ in Galilean space $G_{3}.$ Then it is easy to
obtain that $\tau _{1}=0.$

\begin{center}
\textbf{4. Mannheim Curves in Pseudo - Galilean Space}
\end{center}

\textbf{Definition 4.1. }Let $\alpha $ and $\alpha _{1}$ be an admissible
curves defined in type I with nonzero $\kappa _{\alpha },$ $\tau _{\alpha }$%
, $s\in I$ in Pseudo - Galilean space and $\{T_{\alpha },N_{\alpha
},B_{\alpha }\}$ and $\{T_{\alpha _{1}},N_{\alpha _{1}},B_{\alpha _{1}}\}$
be Frenet frame in Pseudo - Galilean space $G_{3}^{1}$ along $\alpha $ and $%
\alpha _{1},$ respectively. If there exists a corresponding reletionship
between the admissible curves $\alpha $ and $\alpha _{1}$ such that, at the
corresponding points of the admissible curves, principal normal lines $%
N_{\alpha }$ of $\alpha $ coincides with the binormal lines $B_{\alpha _{1}}$
of $\alpha _{1}$, then $\alpha $ is called an admissible Mannheim curves and 
$\alpha _{1}$ is called an admissible Mannheim partner curve of $\alpha .$
The pair $\{\alpha ,\alpha _{1}\}$ is said to be an admissible Mannheim pair
in Pseudo - Galilean space $G_{3}^{1}.$ [1]

\textbf{Theorem 4.2.} Let $\alpha $ an admissible curve defined in Type I in
Pseudo - Galilean space $G_{3}^{1}.$ Then $\alpha $ is an admissible
Mannheim curve if and only if its curvature $\kappa _{\alpha }$ and torsion $%
\tau _{\alpha }$ satisfy the relation $\kappa _{\alpha }=-c\tau _{\alpha
}^{2}$ for some constant $c.$

Proof. Let $\alpha =\alpha (s)$ be an admissible Mannheim curve. Let us
denote by $\{T_{\alpha },$ $N_{\alpha },$ $B_{\alpha }\}$ the Frenet frame
field of $\alpha .$

Assume that $\alpha _{1}=\alpha _{1}(s_{1})$ is an admissible curve whose
binormal \linebreak direction coincides with the principal normal of $\alpha
.$ Namely let us denote by $\{T_{\alpha _{1}},N_{\alpha _{1}},B_{\alpha
_{1}}\}$ the Frenet frame field of $\alpha _{1}.$ Then $B_{\alpha
_{1}}(s_{1})=\pm $ $N_{\alpha }(s).$

The curve $\alpha _{1}$ is parametrized by arclength $s$ as%
\begin{equation}
\alpha _{1}(s)=\alpha (s)+\lambda (s)N_{\alpha }(s)  \tag{4.1}
\end{equation}%
for some function $\lambda (s)\neq 0.$ Differentiating (4.1) with respect to 
$s,$ we find 
\begin{equation}
\alpha _{1}^{^{\prime }}=T_{\alpha }+\lambda ^{\prime }N_{\alpha }+\lambda
\tau _{\alpha }B_{\alpha }.  \tag{4.2}
\end{equation}%
Since the binormal direction of $\alpha _{1}$ coincides with the principal
normal of $\alpha ,$ we have $\lambda ^{\prime }=0.$ Hence $\lambda $ is
constant. The second derivative $\alpha _{1}^{^{^{\prime \prime }}}$ with
respect to $s$ is 
\begin{equation}
\alpha _{1}^{^{^{\prime \prime }}}=(\kappa _{\alpha }+\lambda \tau _{\alpha
}^{2})N_{\alpha }+\lambda \tau _{\alpha }^{\prime }B_{\alpha }.  \tag{4.3}
\end{equation}

Since $N_{\alpha }$ is in the binormal direction of $\alpha _{1},$ we have 
\[
\kappa _{\alpha }+\lambda \tau _{\alpha }^{2}=0. 
\]

Conversely, let $\alpha $ be an admissible curve. Then the curve%
\[
\alpha _{1}(s)=\alpha (s)+\lambda N_{\alpha }(s) 
\]%
has binormal direction $N_{\alpha }.$

\bigskip

\textbf{Theorem 4.3.} Let $\beta $ an admissible curve defined in Type II in
Pseudo - Galilean space $G_{3}^{1}.$ Then $\beta $ is an admissible Mannheim
curve if and only if its curvature $\kappa _{\beta }$ and torsion $\tau
_{\beta }$ satisfy the relation $\kappa _{\beta }=\frac{-c}{\cosh \phi }\tau
_{\beta }^{2}$ for some constant $c.$

Proof. If we\textbf{\ }consider equations (2.3) and proof of the Theorem
4.2. we can prove the theorem easily.

\bigskip \textbf{Theorem 4.4. }Let $\alpha $ be an admissible Mannheim curve
defined by Type I in Pseudo - Galilean space $G_{3}^{1}.$ Then $\alpha _{1}$
is the admissible Mannheim partner curve of $\alpha .$ Then the curvature $%
\kappa _{1}$ and the torsion $\tau _{1}$ of $\alpha _{1}$ satisfy the
following equation

\begin{equation}
\tau _{1}^{^{\prime }}=\frac{\kappa _{1}}{\lambda }(\lambda ^{2}\tau
_{1}^{2}-1)  \tag{4.4}
\end{equation}%
for some nonzero constant $\lambda .$

Proof. It is similar to proof of theorem 3.3.

\textbf{Remark 4.5.} By a simple parameter transformation, the condition 
\[
\tau _{1}^{^{\prime }}=\frac{\kappa _{1}}{\lambda }(\lambda ^{2}\tau
_{1}^{2}-1) 
\]%
can be written as

\[
\tau _{1}=-\frac{\varepsilon }{\lambda }\tan (\varepsilon \dint \kappa
_{1}ds+c_{0}). 
\]%
Therefore, for each an admissible Mannheim curve in Pseudo - Galilean space $%
G_{3}^{1},$ there is an unique Mannheim partner curve.

This reality is true for Mannheim curve in Galilean space $G_{3}.$

\textbf{Proposition 4.6..} Let $\alpha $ be an admissible Mannheim curve in
Pseudo - Galilean space $G_{3}^{1}$ and $\alpha _{1}$ be the admissible
Mannheim partner curve of $\alpha .$ If $\alpha $ is a generalized helix,
then $\alpha _{1}$ is a planar curve.

Proof. Let $T,$ $N,$ $B$ the tangent , principal normal and binormal vector
field of the curve $\alpha $, respectively. From the properties of
generalized helix and the definition of admissible Mannheim curves in Pseudo
- Galilean space $G_{3}^{1},$ we have

\[
\kappa <N,P>=0 
\]%
and 
\[
\kappa <B_{1},P>=0 
\]%
for a constant direction $P$ in Pseudo - Galilean space $G_{3}^{1}.$ Then it
is easy to obtain that $\tau _{1}=0.$

\begin{center}
\textbf{References}
\end{center}

[1] Akyi\u{g}it, M., Azak, A.Z. and Tosun, M., \textit{Admissible Mannheim
Curves in Pseudo-Galilean Space }$G_{3},$ arXiv:1001.2440v2 [math.DG] 21 Jan
2010.

[2] Carmo, M.P., \textit{Differential Geometry of Curves and Surfaces},
Pearson Education, 1976.

[3] Divjak, B., \textit{Curves in Pseudo -- Galilean geometry},
Ann.Univ.Sci. Budapest. E\"{o}tv\"{o}s Sect. Math. 41 (1998) 117-128.

[4] Divjak, B. and Sipus, Z.M., \textit{Special curves on ruled surfaces in
Galilean and pseudo-Galilean spaces}, Acta Math. Hungar.,98(3)
(2003),203-215.

[5] Divjak, B. and Sipus, Z.M., \textit{Minding isometries of ruled surfaces
in pseudo-Galilean space}, Journal of Geometry, 77(2003),35-47.

[6] Ersoy, S., Akyi\u{g}it,\ M. and Tosun, M., \textit{A Note an admissible
Mannheim Curves in Galilean Space }$G_{3},$ arXiv:1003.31110v1 [math.DG] 16
Mar 2010.

[7] K\"{u}lahc\i , M.\ and Erg\"{u}t, M., \textit{Bertrand curves of
AW(k)-type in Lorentzian space}, Nonlinear Anal., 70(4), (2009) 1725-1731.

[8] Liu, H. and Wang, F., \textit{Mannheim partner curves in 3-space},
Journal of Geometry, 88 (2008), 120-126.

[9] Orbay, K. and Kasap, E., \textit{On Mannheim partner curves in E%
${{}^3}$%
} , International Journal of Phsical Sciences, 4(5), (2009) 261-264.

[10] \"{O}\u{g}renmi\c{s}, A.O., \"{O}ztekin, H. and Erg\"{u}t, M., \textit{%
Bertrand Curves in Galilean Space and their Characterizations,} Kragujevac
J. Math. 32 (2009) 139-147.

[11] \"{O}ztekin, H.B. and Erg\"{u}t, M., \textit{Null Mannheim Curves in
the Minkowski 3-space E13}, Turk J Math 34 (2010) 1-8.

[12] R\"{o}schel, O., \textit{Die Geometrie des Galileischen raumes}.
Habilitionsschrift, Leoben, 1984.

[13] Yildirim Yilmaz, M. and Bekta\c{s}, M., \textit{General properties of
Bertrand curves in Riemann-Otsuki space}, Nonlinear Anal., 69(10), (2008)
3225-3231.

\end{document}